\documentclass[11pt]{article}

\usepackage{amssymb}
\usepackage{amsmath}

\usepackage{array}

\usepackage{a4wide}

\usepackage{graphicx}
\usepackage{color}
\graphicspath{{Figures/}}
\input{./psfig.sty}
\gdef\path{Figures}

\def\ITEMMACRO #1 ??? #2 ???{\par\vskip4pt\noindent%
\hangindent=#2em\setbox0\hbox{#1\kern4pt}%
\ifdim\wd0<\hangindent\setbox0\hbox to\hangindent{\hss#1\kern7pt}\fi%
\box0\ignorespaces}
  
\def\Item(#1){\ITEMMACRO {\rm (#1)} ??? 1.8 ???}

\def\bItem{\ITEMMACRO \hss$\bullet$ ??? 1.8 ???}

\let\Bitem=\bItem
\def\BrackItem[#1]{\ITEMMACRO [#1] ??? 1.8 ???}
 
\newtheorem{lemma}{Lemma}
\newtheorem{theorem}{Theorem}

\newtheorem{definition}{Definition}

\newtheorem{problem}{Problem}
\newtheorem{claim}{Claim}
\newtheorem{fact}{Fact}

\def\Claim#1.{\medbreak\ni{\bf Claim~#1.}}
\def\Case#1.{\medbreak\ni{\bf Case~#1.}}
\def\SubCase#1.{\medbreak\ni{\bf Subcase~#1.}}

\def\Proof{\ni{\sl Proof.}\ }

\def\qed{\hfill\fbox{\hbox{}}\medskip}

\def\ni{\noindent}

\def\PP{\mathcal{P}}
\def\TT{\mathcal{T}}

\def\PD{\mathcal{D}}

\def\PL{\mathcal{A}}

\def\PA{\mathcal{A}}

\newcommand{\ov}[1]{\overline{#1}}

\def\bga{\kern2pt{\breve{\kern-2pt\gamma}}}

\begin{document}

% *************************************************************
%                   TITLE PAGE
% *************************************************************

\title{
Intersection Graphs of Pseudosegments:\\
 Chordal Graphs\footnote{Abstracts containing part of this work
appeared in the proceedings of GD'06, LNCS 4372, pp.~208--219
and in the proceedings EuroCG08, pp.~67-70.}
}

\author{
      \normalsize Cornelia Dangelmayr\\
      \small\sf Institut f\"ur Mathematik,\\[-1mm]
      \small\sf Freie Universit\"at Berlin\\[-1mm]
      \small\sf {\tt dangel@math.fu-berlin.de}
\and
      \normalsize Stefan Felsner%%
\footnote{Partially supported by DFG grant FE-340/7-1}\\
      \small\sf Institut f\"ur Mathematik,\\[-1mm]
      \small\sf Technische Universit\"at Berlin.\\[-1mm]
      \small\sf {\tt felsner@math.tu-berlin.de}
\and
      \normalsize William T. Trotter\\
      \small\sf School of Mathematics\\[-1mm]
      \small\sf Georgia Institute of Technology\\[-1mm]
       \small\sf {\tt trotter@math.gatech.edu}
}

\date{}
\maketitle
\vbox{}
\vskip-12mm
\vbox{} 
%%
% *************************************************************

\advance\belowcaptionskip-1pt
\def\vminus{\vbox{}\vskip-25pt\vbox{}}

% *************************************************************
 
\begin{abstract}\advance\leftmargin-2cm
\noindent
We investigate which chordal graphs have a representation as
intersection graphs of pseudosegments. For positive we have a
construction which shows that all chordal graphs that can be
represented as intersection graph of subpaths on a tree are
pseudosegment intersection graphs. We then study the limits of
representability.  We describe a family of intersection graphs of
substars of a star which is not representable as intersection graph of
pseudosegments.  The degree of the substars in this example, however,
has to get large.  A more intricate analysis involving a Ramsey
argument shows that even in the class of intersection graphs of
substars of degree three of a star there are graphs that are not
representable as intersection graph of pseudosegments.

Motivated by representability questions for chordal graphs we 
consider how many combinatorially different $k$-segments, i.e.,
curves crossing $k$ distinct lines, an arrangement of $n$ pseudolines
can  host. We show that for fixed $k$ this number is in $O(n^2)$.
This result is based on a $k$-zone theorem for arrangements of pseudolines
that should be of independent interest.
\end{abstract}

%% ***************************************************************************
\section{Introduction}
%% ***************************************************************************

A family of pseudosegments is understood to be a set of  Jordan
arcs in the Euclidean plane that are pairwise either disjoint or intersect
at a single crossing point. A family of pseudosegments represents a
graph $G$, the vertices of $G$ are the Jordan arcs and two vertices
are adjacent if and only if the corresponding arcs intersect.
A graph represented by a family of pseudosegments is a
{\em pseudosegment intersection graph}, for short a PSI-graph.

PSI-graphs are sandwiched between the larger class of string-graphs
(intersection graphs of Jordan arcs without condition on their intersection behavior)
and of segment-graphs (intersection graphs of straight line segments).
In one of the first papers on the subject Ehrlich et al.~\cite{EET-76} proved
that all planar graphs are string-graphs. In fact this follows from Koebe's
coin graph theorem. Scheinerman in his thesis~\cite{S-84} conjectured
that planar-graphs are segment graphs. Some special cases have been
resolved, most notably de~Castro et al.~\cite{CCDMN-02} proved that
triangle free planar graphs can be represented by segments in three directions.
Recently Chalopin, Gon{\c c}alves and Ochem~\cite{CGO-06} managed to prove that 
all planar graphs are PSI-graphs. 

Interesting research has been conducted regarding the membership
complexity of these classes. The problem for string-graphs was stated
by Graham in~'76.  It was not known to be decidable for almost thirty
years. In two independent papers it was shown that an exponential
number of intersections is sufficient to represent
string-graphs~\cite{PT-02},~\cite{SS-01}.  Later in~\cite{SSS-03} it
was shown that the recognition problem for string-graphs is in NP.
Proofs for NP-hardness have been obtained by Kratochv\ 'il:
In~\cite{K-91} he shows that recognizing string-graphs is NP-hard.
Recognition of PSI-graphs is shown to be NP-complete in~\cite{K-94}.
Recognition of segment-graphs is NP-hard~\cite{KM-89}. Interestingly
it is open whether it is NP-complete. It is known~\cite{KM-94},
however, that a representation via segments with integer endpoints may
require endpoints of size $2^{2^{\sqrt{n}}}$.

There are some classes of graphs where segment representation, hence,
as well PSI-representations, are trivial (e.g.~permutation graphs and
circle graphs) or very easy to find (e.g.~interval graphs).  A large
superclass of interval graphs is the class of chordal graphs.  This
paper originated from investigating chordal graphs in view of their
representability as intersection graphs of pseudosegments.
From the subtree representation of chordal graphs it is immediate
that they are string-graphs.

%% ********************************************************************
\subsection{Basic definitions and results}
%% ***************************************************************************

A graph is {\em chordal} if it has no induced cycles of length greater
than three.  Gavril~\cite{G-74} characterized chordal graphs as the
intersection graphs of subtrees of a tree, that is: a graph
$G=(V_G,E_G)$ is chordal iff there exists a tree $T=(V_T,E_T)$ and a
set $\TT$ of subtrees of $T$ such that there is a mapping $v \to T_v
\in \TT$ with the property that $vw \in E_G$ whenever $T_v \cap T_w
\not =\emptyset$.  The pair $(T,\TT)$ is a {\em tree representation}
of $G$.

A subclass of chordal graphs is the class of
vertex intersection graph of paths on a tree, {\em VPT-graphs} for short.
The precise definition is as follows: A graph $G=(V_G,E_G)$ is a VPT-graph
if there exists a tree $T=(V_T,E_T)$ and a set $\PP$ of paths 
in $T$ such that there is a mapping $v \to P_v \in \PP$ with
the property that $vw \in E_G$ iff $P_v \cap P_w \not =\emptyset$.
Such a pair $(T,\PP)$ is said to be a {\em VPT-representation} of $G$ .

VPT-graphs have been introduced by Gavril in~\cite{G-78}, he called
them path-graph.  Gavril provided a characterization and a recognition
algorithm.  Since then VPT-graphs have been studied continuously,
Monma and Wei~\cite{MW-86} give some applications and many references.
We show:

\begin{theorem}\label{thm:vpt-is-psi}
Every  VPT-graph has a PSI-representation.
\end{theorem} 

The proof of the theorem is given in Section~\ref{sec:vpt-in-psi}.  At
this point we content ourselves with an indication that the result is
not as trivial as it may seem at first glance. Let a
VPT-representation $(T,\PP)$ of a graph $G$ be given. If we fix a
plane embedding of the tree we obtain an embedding of each of the
paths $P_v$ corresponding to $v\in V_G$.  The first idea for
converting a VPT-representation into a PSI-representation would be to
slightly perturb $P_v$ into a pseudosegment~$s_v$ and make sure that
paths with common vertices intersect exactly once and are disjoint
otherwise.  Figure~\ref{fig:intuitive} gives an example of a set of
three paths which can't be perturbed locally as to give a
PSI-representation of the corresponding subgraph.

%%%%%%%%%%%%%%%%%%%%%%%%%%%%%%%%%%%%%%%%%%%%%%%%%%%%%%%%%%%
%%
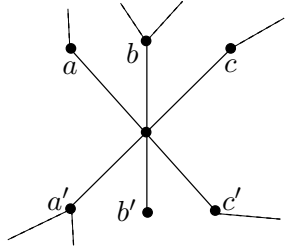
\PsFigCap{25}{intuitive}{The paths $P(a,a')$, $P(b,b')$ and $P(c,c')$ 
can't be perturbed locally to yield a PSI-representation.}
%%
%%%%%%%%%%%%%%%%%%%%%%%%%%%%%%%%%%%%%%%%%%%%%%%%%%%%%%%%%%%

It is natural to ask whether all chordal graphs are PSI-graphs.  This
is not the case, in Theorems~\ref{thm:K-3-33-not-psi} and
\ref{thm:indep+3-not-psi} we give examples of chordal graphs which are
not PSI-representable.  Our examples are quite big, all of them have
more than 5000 vertices.  This is in contrast to the size of
obstructions against PSI-representability within the class of all
graphs. There we know of quite small examples, e.g., the complete
subdivisions of~$K_5$ and~$K_{3,3}$ with 15 vertices.

In Theorem~\ref{thm:K-3-33-not-psi} we show that certain graphs~$K_n^3$
are not in PSI for all $n\geq 33$. The vertex set of~$K_n^3$ is partitioned
as $V=V_C \cup V_I$ such that $V_C=[n]$ induced a clique
and $V_I = {[n]\choose 3}$ is an independent set. The edges between
$V_C$ and $V_I$ represent membership, i.e, $\{i,j,k\} \in V_I$ is
connected to the vertices $i$, $j$ and $k$ from $V_C$.

The graph~$K_n^3$ can be represented as intersection graph of subtrees
of a star~$S$ with ${n \choose 3}$ leaves. The leaves correspond to
the elements of $V_I$. Each $v\in V_I$ is represented by trivial tree
with only one node. A vertex $i$ of the complete graph is represented
by the star connecting to all leaves of triples containing $i$. This
representation shows that~$K_n^3$ is chordal.  The central node of the
star $S$ has high degree. If we take a path of ${n \choose 3}$ nodes
and attach a leaf-node to each node of the path we obtain a tree $T$
of maximum degree three such that the graph~$K_n^3$ can be represented
as intersection graph of subtrees of $T$.  Actually the tree $T$ and
its subtrees are caterpillars of maximum degree three.

These remarks show that the positive result of Theorem~\ref{thm:vpt-is-psi} and
the negative of Theorem~\ref{thm:K-3-33-not-psi} only leave a small
gap for questions:
If the subtrees, in a tree representation of a chordal graph, are paths 
we have a PSI-representation. If we allow the subtrees to be
stars (of large degree), or (large) caterpillars of maximum degree three,
there need not exist a PSI-representation. 

What if we restrict the host tree to be a star and the subtrees to be
substars of degree three? Let $S_n$ be the chordal graph whose
vertices are represented by all substars with three leaves and all
leaves on a star with $n$ leaves.  In
Theorem~\ref{thm:indep+3-not-psi} we show that~$S_n$ is not
PSI-representable if $n$ is large enough. This resolves our major
conjecture from~\cite{DF-06}. The proof makes use of a Ramsey
argument, hence, we need a really large $n$ for the result. Note that
as in the case of~$K_n^3$ the vertex set of $S_n$ again partitions
into a clique $V_C$ and an independent set $V_I$.

In Section~\ref{sec:k-seg} we deal with graphs whose vertex set splits
into a clique $V_C$ of size $n$ and a set $W$. We impose no condition
on the subgraph induced by $W$ but require that each~$w\in W$ has
three neighbors in $V_C$ and no two vertices in $W$ have the same
neighbors in $V_C$. We conjecture that if such a graph has a
PSI-representation then $|W| \in o(n^{2-\epsilon})$. Motivated by this
problem we consider how many combinatorially different $k$-segments,
i.e., curves crossing~$k$ distinct lines, an arrangement of $n$
pseudolines can host. In Theorem~\ref{thm:k-seg} we show that for
fixed $k$ this number is in $O(n^2)$.  This result is based on a
$k$-zone theorem for arrangements of pseudolines
(Theorem~\ref{thm:k-zone}).  This theorem states that the complexity
of the $k$-zone in an arrangement of $n$ pseudolines is in $O(n)$.

We conclude with some open problems.

%  ***************************************************************
\section{A PSI-representation for VPT-graphs}
\label{sec:vpt-in-psi}
%% ***************************************************************************

\subsection{Preliminaries}

A path in a tree $T$ with endpoints $a$ and $b$ is denoted $P(a,b)$. 
If both endpoints of a path $P$ in $T$ are leaves
we call $P$ a {\it leaf-path}.

\begin{lemma} 
\label{L1}
Every VPT-graph has a path-representation $(T,\PP)$ such that all paths in
$\PP$ are leaf-paths and no two vertices are represented by the same leaf-path. 
\end{lemma}

\Proof
 Let an arbitrary VPT-representation $(T,\PP)$ of
 $G$ with $P_v=P(a_v,b_v)$ for all $v \in V_G$ be given. 
 Now let $\ov{T}$ be the tree obtained from $T$ by attaching a new
 node $\ov{x}$ to every node $x$ of $T$. Representing the vertex
 $v$ by the path $\ov{P_v}=P(\ov{a_v},\ov{b_v})$ in $\ov{T}$ yields
 a VPT-representation of $G$ using only leaf-paths.
\qed

As a class of intersection graphs the class of VPT-graphs is closed under taking
induced subgraphs.This observation together with the previous lemma
show that Theorem~\ref{thm:vpt-is-psi}
is implied by the following:

\begin{theorem}
\label{t3}
Given a tree $T$ we let
$G$ be the VPT-graph whose vertices are in bijection to the set
of all leaf-paths of $T$. The graph $G$ has a PSI-representation
with pseudosegments $s_{i,j}$ corresponding to the paths
$P_{i,j}=P(l_i,l_j)$ in $T$. In addition there is
a collection of pairwise disjoint disks, 
one disk $R_i$ associated with each leaf $l_i$ of $T$,
such that:
\Item(a) The intersection $s_{i,j}\cap R_k \not=\emptyset$ if and only if $k=i$ or $k=j$. Furthermore the intersections $s_{i,j}\cap R_i$ and $s_{i,j}\cap R_j$ are Jordan curves.
\Item(b) Any two pseudosegments intersecting $R_i$ 
cross in the interior of this disk.
\end{theorem}

\ni
We will prove Theorem~\ref{t3} by induction on the number of inner
nodes of tree~$T$. The construction will have multiple intersections,
i.e., there are points where more than two pseudosegments intersect.
By perturbing the pseudosegments participating in a multiple
intersection locally the representation can easily be transformed into
a representation without multiple intersections.

%%%%%%%%%%%%%%%%%%%%%%%%%%%%%%%%%%%%%%%%%%%%%%%%%%%%%%%%%%
\subsection{Theorem~\ref{t3} for trees with
  one inner node}\label{sec:IndAnf}

Let $T$ have one inner node $v$ and let $L=\{l_1,..,l_m\}$ be the set
of leaves of $T$.  The subgraph~$H$ of $G$ induced by the set
$\PP=\{P(l_i,l_j)\mid l_i,l_j \in L, l_i\not=l_j\}$ of leaf-paths is a
complete graph on~$m \choose 2$ vertices, this is because every path
in $\PP$ contains $v$.

Take a circle $\gamma$ and choose $m$ points $c_1,..,c_m$ on $\gamma$
such that the set of straight lines spanned by pairs of different
points from $c_1,..,c_m$ contains no parallel lines.
For each $i$ choose a small disk $R_i$ centered at $c_i$ 
such that these disks are disjoint and put them in one-to-one
correspondence with the leaves of $T$.
Let $s_{i,j}$ be the line connecting $c_i$ and $c_j$. If the disks
$R_k$ are small enough we clearly have :
\Item(a) 
The line $s_{i,j}$  intersects  $R_i$ and $R_j$ but no further disk $R_k$.
\Item(b) 
Two lines $s_I$ and $s_J$ with $I\cap J = \{i\}$ contain corner $c_i$, hence $s_I$ and $s_J$
cross in the disk~$R_i$.
\medskip

%%%%%%%%%%%%%%%%%%%%%%%%%%%%%%%%%%%%%%%%%%%%%%%%%%%%%%%%%%%
%%
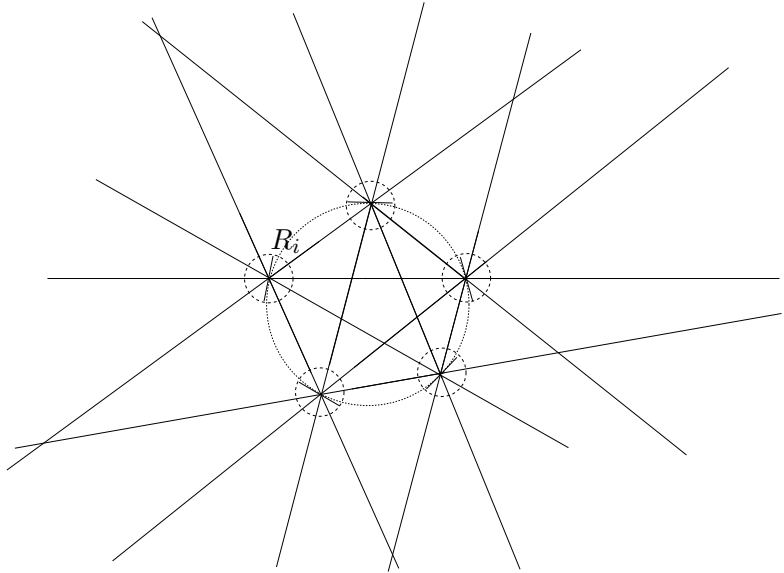
\PsFigCap{15}{circle-two}{The construction for the star with five leaves.}
%%
%%%%%%%%%%%%%%%%%%%%%%%%%%%%%%%%%%%%%%%%%%%%%%%%%%%%%%%%%%%

\ni
Prune the lines such that the remaining part of each $s_{i,j}$
still contains its intersection with all the other lines and all
segments have there endpoints on a circumscribing circle $C$. Every
pair of segments stays intersecting, hence, we have a segment intersection
representation of~$H$.

Add a diameter $s_{i,i}$ to every disk $R_i$, this segment serves as
representation for the leaf-path~$P_{i,i}$. Altogether we have
constructed a representation of $G$ obeying the required properties
(a) and (b), see Figure~\ref{fig:circle-two} for an illustration.

%% ********************************************************************
\subsection{Theorem~\ref{t3} for trees with more than one inner node}

Now let $T$ be a tree with inner nodes
$N=\{v_1,..,v_n\}$ and assume the
theorem has been proven for trees with at most $n-1$ inner
nodes. Let $L=\{ l_1,..,l_m\}$ be the set of
leaves of $T$. With $L_i\subset L$ we denote the set of leaves attached to $v_i$.
We have to produce a PSI-representation of the intersection graph $G$
of  $\PP=\{P_{i,j} \mid l_i,l_j \in L \}$, i.e., of the set of all
leaf-paths of $T$. We suppose that $v_n$ is a leaf node in the tree induced by $N$: 

\Bitem 
The tree $T_n$ is 
the star with inner node $v_n$ and its leaves $L_n=\{l_k,..,l_m\}$.

\Bitem
The tree $T'$ contains all
nodes of $T$ except the leaves in $L_n$. The
set of inner nodes of $T'$ is $N'=N\backslash \{ v_n\}$, the
set of leaves is $L'=L \backslash L_n \cup \{v_n\}$.
For consistency we rename $l_0:=v_n$ in $T'$,
hence $L'=\{l_0,l_1,..,l_{k-1}\}$.
\medskip

%%%%%%%%%%%%%%%%%%%%%%%%%%%%%%%%%%%%%%%%%%%%%%%%%%%%%%%%%%%
%%
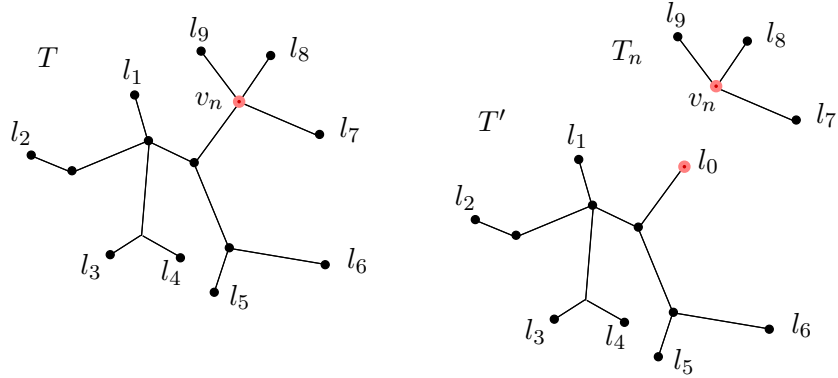
\PsFigCap{25}{tree-two}{A tree $T$ and the two induced subtrees $T'$
  and $T_n$.}
%%
%%%%%%%%%%%%%%%%%%%%%%%%%%%%%%%%%%%%%%%%%%%%%%%%%%%%%%%%%%%

\ni
Let $G_n$ and $G'$ be the VPT-graphs induced by all leaf-paths in
$T_n$ and $T'$.  Both these trees have fewer inner nodes than
$T$. Therefore, by induction we can assume that we have
PSI-representations $PS'$ of $G'$ and $PS_n$ of $G_n$ as claimed in
Theorem~\ref{t3}.
We will construct a PSI-representation of $G$ using $PS'$ and $PS_n$.
The idea is as follows:

\begin{enumerate}
\item
  Replace every pseudosegment of $PS'$ representing a leaf-path
  ending in $l_0$ by a bundle of pseudosegments. This bundle stays
  within a narrow tube around the original pseudosegment.
\item
  Remove all pieces of pseudosegments
  from the interior of the disk $R_0$ and
  patch an appropriately transformed copy of $PS_n$ into $R_0$.
\item
  The crucial step is to connect the pseudosegments of the
  bundles through the interior of $R_0$ such that the induction invariants
  for the transformed disks of $R_r$ with $k\leq r\leq m$ are satisfied. 
\end{enumerate}

\ni
The set $\PP$ of leaf-paths of $T$ can be partitioned into three
parts.  The subsets~$\PP'$ and $\PP_n$ are leaf-paths of $T'$ or $T_n$
let the remaining subset be~$\PP^*$. The paths in~$\PP^*$ connect
leaves $l_i$ and~$l_r$ with $1\leq i < k \leq r \leq m$, in other
words they connect a leaf $l_i$ from $T'$ through $v_n$ with a leaf in
$T_n$. We subdivide these paths into classes $\PP_1^*,..,\PP_{k-1}^*$
such that $\PP_i^*$ consists of those paths from~$\PP^*$ which start
in the leaf $l_i$ of $T'$.  Each $\PP_i^*$ consists of $|L_n|$ paths. In $T'$ we have
the pseudosegment~$s_{i,0}$ which leads from $l_i$ to $l_0$. Replace
each such pseudosegment~$s_{i,0}$ by a bundle of $|L_n|$ parallel pseudosegments
routed in a narrow tube around~$s_{i,0}$.

We come to the second step of the construction.  Remove all pieces of
pseudosegments from the interior of $R_0$. Recall that the
representation $PS_n$ of $G_n$ from~\ref{sec:IndAnf} has the property
that all long pseudosegments have their endpoints on a circle $C$.
Choose two arcs $A_b$ and $A_t$ on $C$ such that every segment spanned
by a point in $A_b$ and a point in $A_t$ intersects each pseudosegment
$s_{i,j}$ with $i\neq j$, this is possible by the choice of $C$. This
partitions the circle into four arcs which will be called
$A_b,A_l,A_t,A_r$ in clockwise order. The choice of $A_b$ and $A_t$
implies that each pseudosegment touching $C$ has one endpoint in $A_l$
and the other in~$A_r$.

Map the interior of $C$ with an homeomorphism $h$ into a wide
rectangular box $\Gamma$ such that~$A_t$ and~$A_b$ are mapped to the
top and bottom sides of the box, $A_l$ is the left side and~$A_r$ the
right side. This makes the images of all long pseudosegments traverse
the box from left to right. We may also require that the homomorphism
maps the disks $R_r$ to disks and arranges them in a nice left to
right order in the box, Figure~\ref{fig:box} shows an example. The
figure was generated by sweeping the representation from
Figure~\ref{fig:circle-two} and converting the sweep into a wiring
diagram (the diametrical segments $s_{r,r}$ have been re-attached
horizontally).

%%%%%%%%%%%%%%%%%%%%%%%%%%%%%%%%%%%%%%%%%%%%%%%%%%%%%%%%%%%
%%
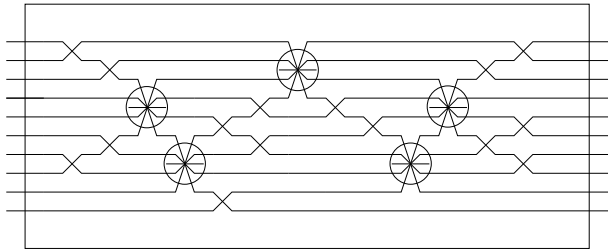
\PsFigCap{25}{box}{A box containing a deformed copy of the
  representation from Figure~\ref{fig:circle-two}. }
%%
%%%%%%%%%%%%%%%%%%%%%%%%%%%%%%%%%%%%%%%%%%%%%%%%%%%%%%%%%%%

In the box we have a left to right order of the disks
$R_r$, $l_r \in L_n$. By possibly relabeling the leaves of $L_n$ we can
assume that the disks are ordered from left to right as
$R_k,..,R_m$. 

This step of the construction is completed by placing the box $\Gamma$
appropriately resized in the disk $R_0$ such that each of the
segments~$s_{i,0}$ from the representation of $G'$ traverses the box
from bottom to top and the sides $A_l$ and $A_r$ are mapped
to the boundary of $R_0$. The boundary of $R_0$ is thus partitioned
into four arcs which are called $A_l,A'_t,A_r,A'_b$ in clockwise order.
We assume that the segments~$s_{i,0}$ touch the
arc $A'_b$ in $PS'$ in counterclockwise order as $s_{1,0},..,s_{k-1,0}$,
this can be achieved by renaming the leaves appropriately. 

%%%%%%%%%%%%%%%%%%%%%%%%%%%%%%%%%%%%%%%%%%%%%%%%%%%%%%%%%%%%%%%%%%%%%%%%%%%%

Note that by removing everything from the interior of $R_0$ we have
disconnected all the pseudosegments which have been inserted in
bundles replacing the original pseudosegments~$s_{i,0}$. Let
$B_{i}^{in}$ be the half of the bundle of~$s_{i,0}$ which
touches~$A'_b$ and let $B_{i}^{out}$ be the half which touches $A'_t$.
By the above assumption the bundles $B_{1}^{in},..,B_{k-1}^{in}$ touch
$A'_b$ in counterclockwise order, with (b) from the statement of the
theorem and induction it follows that $B_{1}^{out},..,B_{k-1}^{out}$
touch~$A'_t$ in counterclockwise order. Within a bundle $B_{i}^{in}$
we label the segments as $s_{i,k}^{in},..,s_{i,m}^{in}$, again
counterclockwise. The segment in $B_{i}^{out}$ which was connected
to~$s_{i,r}^{in}$ is labeled $s_{i,r}^{out}$ The pieces $s_{i,r}^{in}$
and $s_{i,r}^{out}$ will be part of the pseudosegment representing the
path~$P_{i,r}$.

To have property (b) for the pseudosegments of a bundle we twist
whichever of the bundles~$B_i^{in}$ or~$B_i^{out}$ traverses $R_i$
within this disk $R_i$ thus creating a multiple intersection point.
Note that they all cross $s_{i,i}$ as did $s_{i,0}$.

Also due to (b) the pseudosegments of paths $P_{i,r}$ for fixed $r \in
\{k,..,m\}$ have to intersect in the disks $R_r$ inside of the
box $\Gamma$. To prepare for this we take a narrow bundle of $k-1$ parallel
vertical segments reaching from top to bottom of the box $\Gamma$ and
intersecting the disk~$R_r$.  This bundle is twisted in the
interior of $R_r$. Let $\check{a}^r_1,..,\check{a}^r_{k-1}$ be
the bottom endpoints of this bundle from left to right and let
$\hat{a}^r_1,..,\hat{a}^r_{k-1}$ be the top endpoints from right to
left, due to the twist the endpoints $\check{a}^r_j$ and $\hat{a}^r_j$
belong to the same pseudosegment.

We are ready now to construct the pseudosegment
$s_{i,r}$ that will represent the path $P_{i,r}$ in $T$ for $1\leq i < k \leq r \leq m$. The first part of $s_{i,r}$
is $s_{i,r}^{in}$, this pseudosegment is part of the bundle
$B_i^{in}$ and has an endpoint on $A'_b$. Connect this
endpoint with a straight segment to $\check{a}^r_i$,
from this point there is the connection up to $\hat{a}^r_i$.
This point is again connected by a straight segment to the endpoint of $s_{i,r}^{out}$
on the arc $A_t'$. The last part of $s_{i,r}$ is the pseudosegment
$s_{i,r}^{out}$ in the bundle $B_{i}^{out}$.
The construction is illustrated in Figure~\ref{fig:rewire2}.

%%%%%%%%%%%%%%%%%%%%%%%%%%%%%%%%%%%%%%%%%%%%%%%%%%%%%%%%%%%
%%
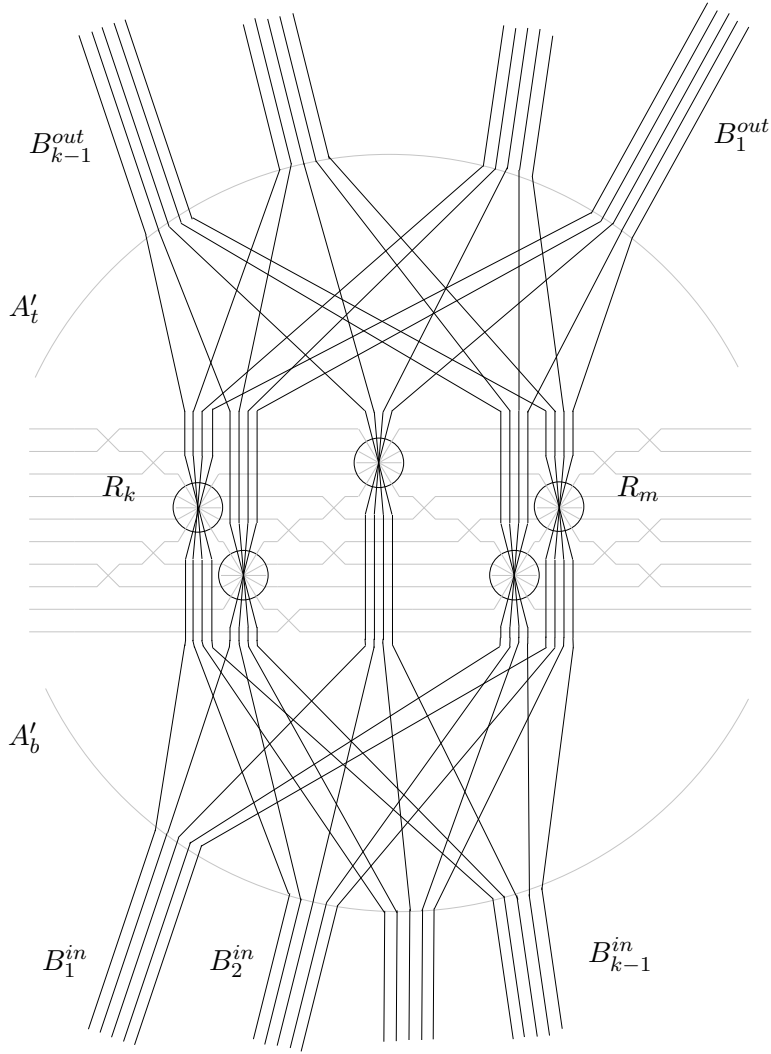
\PsFigCap{30}{rewire2}{The routing of pseudosegments in the disk $R_0$, an example.}
%%
%%%%%%%%%%%%%%%%%%%%%%%%%%%%%%%%%%%%%%%%%%%%%%%%%%%%%%%%%%%

The following list of claims collects some crucial properties of 
the construction.

\Claim 1. There is exactly one pseudosegment $s_{i,j}$ for every pair
$l_i,l_j$ of leaves of $T$.

\Claim 2. The pseudosegment  $s_{i,j}$ traverses  $R_i$ and $R_j$ but
stays disjoint from every other of the disks.

\Claim 3. Any two pseudosegments intersecting the disk $R_i$ cross in
$R_i$. 
  
\Claim 4. Two pseudosegments $s_{i,j}$ and $s_{i',j'}$ intersect
exactly if the corresponding paths $P_{i,j}$ and $P_{i',j'}$ intersect
in $T$.

\Claim 5. Two pseudosegments $s_{i,j}$ and $s_{i',j'}$ intersect at
most once, i.e., to call them pseudosegments is justified.
 
\medskip

\ni
These claims follow by induction.
Hence the construction indeed yields a
representation of $G$ as intersection graph of pseudosegments
and this representation has properties (a) and (b).

%% ********************************************************************
%% ********************************************************************
\section{Chordal graph that are not PSI}
\label{sec:proof2}
%% ***************************************************************************
%% ********************************************************************

Recall the definition of the graphs~$K_n^3$:
The graph~$K_n^3$ has two groups of vertices. The set $V_C = [n]$ induces a
clique in addition every triple $\{i,j,k\} \subset
[n]$ is a vertex adjacent only to the three 
vertices $i$, $j$ and $k$, hence, the triples form an independent set $V_I$.
In the introduction we have already described tree-representations of~$K_n^3$,
in particular~$K_n^3$ is chordal.

\begin{theorem}\label{thm:K-3-33-not-psi}
For $n\geq 33$ the graph~$K_n^3$ admits no PSI-representation.
\end{theorem}

Assuming that there is a representation of~$K_n^3$ as intersection
graph of pseudosegments. Let $PS_C$ and $PS_I$ be the sets of
pseudosegments representing vertices from $V_C$ and $V_I$.

The pseudosegments in $PS_C$ form a set of pairwise crossing
pseudosegments, we refer to the configuration of these pseudosegments as the 
arrangement $A_n$. The set $S = PS_I$ of {\em small} pseudosegments has the following
properties:  

\Item(i) 
Any two pseudosegments $t\neq t'$ from $S$ are disjoint,
\Item(ii)
Every pseudosegment $t\in S$ has nonempty intersection with exactly three
pseudosegments from the arrangement $A_n$ and no two pseudosegments
$t\neq t'$ intersect the same three pseudosegments from $A_n$.
\medskip

\ni
The idea for the proof is to show that a set of pseudosegments with
properties (i) and (ii) only has $O(n^2)$ elements. The theorem
follows, since $|S| = {n \choose 3} = \Omega(n^3)$.

%% ********************************************************************
\subsection{$K_n^3$ and planar graphs} 
%% ********************************************************************

Every pseudosegment $p \in A_n$ is cut into $n$ pieces
by the $n-1$ other pseudosegments of $A_n$. Let $W$ be the set
of all the pieces obtained from pseudosegments from~$A_n$,
note that $|W|=n^2$.  Pseudosegments in $S$ intersect
exactly three pieces from three different pseudosegments of $A_n$. 
Elements of $S$ are called {\em 3-segments}. Every 3-segment has a unique middle and
two outer intersections. Let $S(w)$ be the
set of 3-segments with middle intersection on the piece~$w \in W$.
This yields the partition $S=\bigcup_{w\in W}S(w)$ of the set of 3-segments.

Define $G_p=(W,E_p)$ as the simple graph where two pieces $w,w'$ are
adjacent if and only if there exists a 3-segment $t\in S$
such that $t$ has its middle intersection on
$w$ and an outer intersection on $w'$.

\begin{lemma}
 $G_p=(W,E_p)$ is planar. 
\end{lemma}  
\Proof
A planar embedding of $G_p$ is induced by $A_n$ and $S$.  Contract all
pieces from $A_n$, the contracted pieces represent
the vertices of $G_p$. The 3-segments are pairwise
non-crossing, this property is maintained during contraction of
pieces, see Figure~\ref{fig:Gp-neu}. If a 3-segment $t\in S$ has middle piece $w$
and outer pieces $w'$ and~$w''$, then $t$ contributes the two edges
$(w,w')$ and $(w,w'')$ to $G_p$. Hence, the multi-graph obtained through
these contractions is planar and its underlying simple graph is indeed $G_p$.
\qed

%%%%%%%%%%%%%%%%%%%%%%%%%%%%%%%%%%%%%%%%%%%%%%%%%%%%%%%%%%%
%%
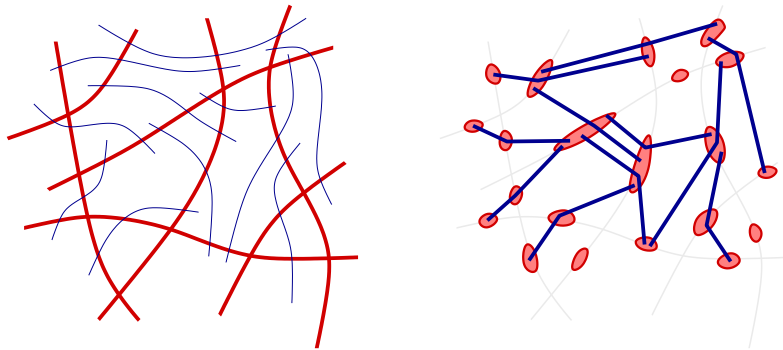
\PsFigCap{30}{Gp-neu}{A part of $A_n$ with some 3-segments and
                      the edges of $G_p$ induced by them. }
%%
%%%%%%%%%%%%%%%%%%%%%%%%%%%%%%%%%%%%%%%%%%%%%%%%%%%%%%%%%%%

Let $N(w)$ be the set of neighbors of $w \in W$ in $G_p$ and let $d_{G_p}(w)=|N(w)|$.

\begin{lemma}
  \label{deg} 
  The size of a set $S(w)$ of 3-segments with
  middle piece $w$ is bounded by $3 d_{G_p}(w)-6$ for every $w\in
  W$.\end{lemma}
 
\Proof The idea is to contract just the pieces
corresponding to elements of $N(w)$ to points. The 3-segments
in $S(w)$ together with the vertices obtained by contraction form a
planar graph. Note that this time the graph is not a multi-graph,
since every 3-segment leading to an edge also intersects the
piece $w$. The number of vertices of this graph is $d_{G_p}(w)$,
hence, there are at most $3 d_{G_p}(w)-6$ edges. \qed

In~\cite{DF-06} we have shown the stronger statement, that the graph $N_w$
is acyclic, hence a forest. This implies that the number of edges is
at most $d_{G_p}(w)-1$. We omit the more involved proof. Still we
will use the stronger bound in the following calculations.
Following the computation given below the weaker bound of Lemma~\ref{deg}
implies that~$K_n^3$ is not PSI-representable for $n \geq 75$.

%%%%%%%%%%%%%%%%%%%%%%%%%%%%%%%%%%%%%%%%%%%%%%%%%%%%%%%%%%%
%%
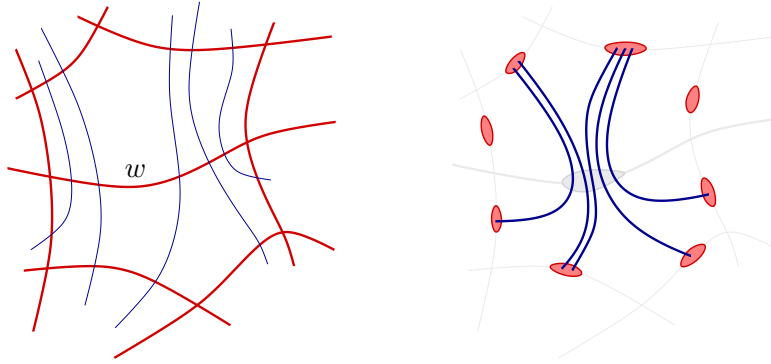
\PsFigCap{19}{vtb-neu}{The planar graph $NG_w$ induced by the
  3-segments with middle intersection. }
%%
%%%%%%%%%%%%%%%%%%%%%%%%%%%%%%%%%%%%%%%%%%%%%%%%%%%%%%%%%%%

\ni
Recall that from simple counting and from the planarity of $G_p$ we have:
\Bitem 
$|W| = n^2$,
\Bitem 
$\sum_{ w \in W}   d_{G_p}(w) = 2 | E_p| < 6|V_p| = 6|W|$.
\bigskip

\ni
Using $|S(w) \leq d_{G_p}(w)-1$ we thus obtain:
\Bitem 
$|S| =\sum_{w \in W} |S(w)| \leq \sum_{w \in W} (d_{G_p}(w)-1) 
< 6 |W| -|W|  = 5n^2$.
\bigskip

\ni
Since ${n \choose 3} >5n^2$ for all $n \geq 33$
we conclude that~$K_n^3$ does not belong to PSI for $n \geq 33$.
This completes the proof of Theorem~\ref{thm:K-3-33-not-psi}. 
\qed

%%%%%%%%%%%%%%%%%%%%%%%%%%%%%%%%%%%%%%%%%%%%%%%%%%%%%%%%%%%%%%%%%%%%
\section{$k$-Segments in Arrangements of Pseudolines}
\label{sec:k-seg}
%%%%%%%%%%%%%%%%%%%%%%%%%%%%%%%%%%%%%%%%%%%%%%%%%%%%%%%%%%%%%%%%%%%%

In the previous section we have shown that the chordal graph $K_n^3$,
$n\geq 33$, has no PSI-representation. For this purpose we partitioned
the graph $K_n^3$ into two induced subgraphs, the maximal clique and
the maximal independent set. In any PSI-representation of $K_n^3$, the
clique has to correspond to an arrangement, that is a set of $n$
pairwise intersecting pseudosegments. Let $\PA_n$ denote a
representation of the clique.  To obtain a PSI-representation of the
whole graph $K_n^3$, it is necessary to assign the $n \choose 3$
vertices of the independent set to a set of disjoint pseudosegments
({\em triple-segments}) in the base arrangement $\PA_n$. 
What if we omit the condition of disjointness for the triple-segments?

To be more precise: 
Let $\PA_n$ be an arrangement of $n$ pseudosegments and 
$T\subset {[n]\choose 3}$ be a set of triples. We say that
{\em $T$ is hosted on $\PA_n$} if for every 
$t \in T\subset {[n]\choose 3}$ there is an $S_t$ such
that $\PA_n \cup \{ S_t : t\in T\}$ is a family of pseudosegments and
the segment $S_t$ with $t=\{i,j,k\}$ is intersecting the lines labeled
$i,j,k$ in $\PA_n$ and no others.

\begin{problem}\label{prob:PS}
What is the growth of the largest function
$f_3(n)$ such that there is a $T$ hosted on some $\PA_n$ with
$f_3(n) = |T|$?
\end{problem}

Although we conjecture that $f_3(n)$ is $o(n^{2+\epsilon})$ for all 
$\epsilon > 0$ we have not even been able to show that 
$f_3(n)$ is $o(n^{3})$. 

The results of this section were motivated by 
Problem~\ref{prob:PS}. We have simplified the situation by
dealing with arrangements of pseudolines instead of arrangements
of pseudosegments. The advantage of this model lies in the fact
that the elements of the arrangement cannot be bypassed at their ends.
We think that this more geometric model is of independent interest.
Indeed our result about the complexity of the $k$-zone 
(Theorem~\ref{thm:k-zone})
has already found applications in the work 
of Scharf~\cite{Sch-08}.

\begin{definition}  
A {\em pseudoline}\/ in the Euclidean plane is a simple curve that
approaches a point at infinity in either direction. An
{\em arrangement of pseudolines} is a family of pseudolines with the
property that each pair of pseudolines has a unique point of
intersection, where the two pseudolines cross. 
\end{definition}

We now ask:
How many triple-segments can an arrangement
$\PL_n$ of $n$ pseudolines host?
\medskip

A pseudoline in an arrangement of pseudolines is split into a sequence
of edges and vertices by the crossing lines.

\begin{definition} 
A {\em $k$-segment} in an arrangement $\PL$ of pseudolines is
a sequence $e_1,e_2,\ldots,e_k$ of edges from $k$ different
pseudolines of $\PL_n$ with the
property that there exists a curve crossing these edges in the given
order that has no further intersections with $\PL$.
\end{definition}

%%%%%%%%%%%%%%%%%%%%%%%%%%%%%%%%%%%%%%%%%%%%%%%%%%%%%%%%%%%
%%
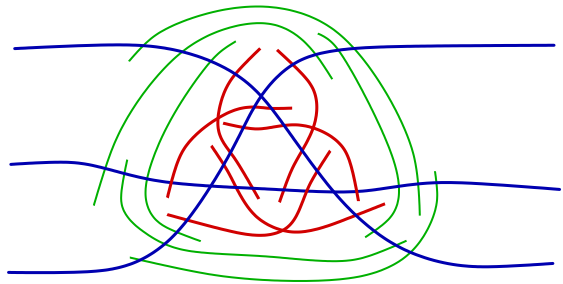
\PsFigCap{40}{arr3+12}{An arrangement of 3 pseudolines and 
              representative curves for all its 3-segments.}
%%
%%%%%%%%%%%%%%%%%%%%%%%%%%%%%%%%%%%%%%%%%%%%%%%%%%%%%%%%%%%

The main theorem of this section is:

\begin{theorem}\label{thm:k-seg}
For $k$ fixed the number of different $k$-segments in an arrangement of 
$n$ pseudolines is at most $c^k n^2 \in O(n^2)$ for some $c$.
\end{theorem}

Compared to the setting of Problem~\ref{prob:PS} we now
consider $k$-segments instead of triple-segments, moreover, we can
count the same set of $k$ pseudolines with many $k$-segments and we
have dropped the condition that the $k$-segments must be compatible,
their representing curves may cross any number of times. 

Note that if we drop the compatibility
restriction in the pseudosegment case all
triples become representable, i.e., the number of representable triples
becomes $\Theta(n^3)$.

The proof of Theorem~\ref{thm:k-seg} will be given in Subsection~\ref{ssec:k-seg}.
It is based on Theorem~\ref{thm:k-zone} which is the main result of
the next subsection.

%%%%%%%%%%%%%%%%%%%%%%%%%%%%%%%%%%%%%%%%%%%%%%%%%%%%%%
\subsection{$k$-zones in arrangements}
%%%%%%%%
\def\EZone{\text{\rm Zone}^+}
\def\FZone{\text{\rm Zone}^-}

\begin{definition} 
  Let $\PL$ be an arrangement of pseudolines, $p\in \PL$ and
  $k\geq 2$. The { \em $k$-zone} of $p$  is defined
  as the collection of vertices, edges and faces of $\PL$ that
  can be connected to $p$ by a curve that has at most 
  $k$ intersections with lines of $\PL\setminus p$.
\end{definition}

We are interested in upper bounds for the complexity of the $k$-zone.
Assume that the pseudoline $p$ is a horizontal straight-line.  Hence,
$p$ splits the plane in two half-planes $H^+$ and~$H^-$ both
containing `one half' of the $k$-zone. More precisely, the double of
an upper bound for the complexity of the intersection of the $k$-zone
with $H^+$ that is valid for all arrangements~$\PL_n$
is an upper bound for the complexity of the $k$-zone.

We are mainly interested in the edges of the $k$-zone.
They can also be defined as the set of edges that
are contained in a $(k+1)$-segment whose initial edge is
on $p$. 
With $\EZone_k(p)$ we denote the set of edges
of the $k$-zone of $p$ that are contained in $H^+$. 

\begin{theorem}\label{thm:k-zone}
The number of edges in  $\EZone_k(p)$ is linear in
$n$, more precisely
$$
|\EZone_k(p)| < 5\cdot 4^{k-1}(n-1)
$$
\end{theorem}

The reminder of this subsection is devoted to the proof of the theorem.
We argue by induction on $k$. The base case $k=1$ is 
just the classical 2-D Zone Theorem from computational
geometry, see e.g.~\cite[page 147]{m-ldg-02}. 

For the induction step we consider 
$Z_k(p) = \EZone_k(p) \setminus \EZone_{k-1}(p)$,
i.e., the set of edges in $H^+$ that can only be connected
to an edge of $p$ via a $k+1$-segment if they
are at one end of the segment and the edge of $p$
is the other end.

\begin{claim}\label{cl:Z_k(p)}
Let $\PL_n$ be an arrangement of
pseudolines and $p\in \PL_n$, then 
$|Z_k(p)| < 14\cdot 4^{k-2}(n-1)$
\end{claim}

This claim and $|\EZone_k(p)| = \sum_{j\leq k} |Z_j(p)|
 <  14(n-1) \sum_{j\leq k} 4^{j-2}$ imply the bound 
stated in the theorem.

The half-space $H^+$
contains a half-line $q^+$ of each pseudoline $q\neq p$ from $\PL_n$.
We orient the half-lines away from $p$, i.e., they become rays
emanating from a point on $p$.

Consider an edge $e \in H^+$ that has no vertex on $p$.  Edge $e$
inherits an orientation from the supporting ray $q^+$. The tail of $e$
is at a crossing of $q^+$ with some other half-line $r^+$. With
respect to the orientation of $r^+$ we can classify $e$
as either {\em right-outgoing} or {\em left-outgoing}.

Let $Z^{r}_k(p)$ be the set of edges  in $Z_k(p)$ that are 
right-outgoing. Any upper bound on the number of right-outgoing
edges in $Z_k(p)$ is also valid for left-outgoing edges.
Hence, the following lemma implies
Claim~\ref{cl:Z_k(p)}

\begin{lemma}
$|Z^{r}_k(p)| < 7\cdot 4^{k-2}(n-1)$
\end{lemma}

\Proof
Consider the edges that are right-outgoing from $q^+$ 
and order them by increasing distance from $p$ as
$e_q^1,e_q^2,\ldots$. 

Although $e_q^k$ need not be in $Z^{r}_k(p)$
we account for the edge $e_q^k$ in our upper bound on
$|Z^{r}_k(p)|$. Summing over all lines this contributes $n-1$.  
The advantage is that
the first $k$ of the outgoing edges from $q^+$ need no further
consideration.

Suppose that there is an $l > k$ such that $e_q^l$ is in $Z^{r}_k(p)$.
A curve that connects $e_q^l$ with $p$ and does not cross~$q^+$ has to
intersect with each of the supporting lines of $e_q^1,\ldots,e_q^k$.
By assumption there is a curve~$\gamma$ connecting $p$ with $e_q^l$
that has exactly $k-1$ crossings with other pseudolines.
Hence,~$\gamma$ has to intersect $q^+$ in some edge $e$, see
Figure~\ref{fig:high-edge}.

%%%%%%%%%%%%%%%%%%%%%%%%%%%%%%%%%%%%%%%%%%%%%%%%%%%%%%%%%%%
%%
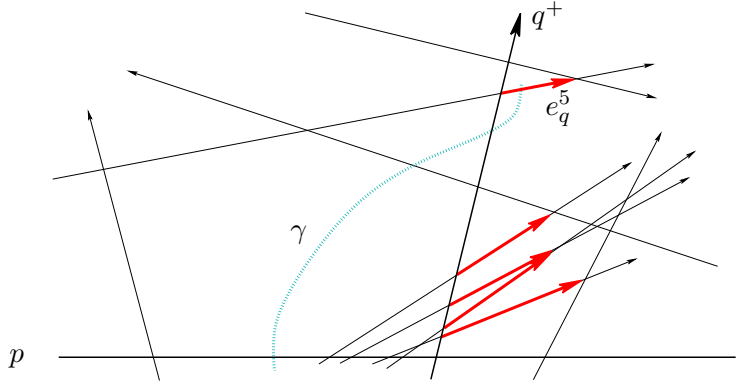
\PsFigCap{30}{high-edge}{The fifth right-outgoing edge of $q^+$
                         is in $Z^{r}_3(p)$.}
%%
%%%%%%%%%%%%%%%%%%%%%%%%%%%%%%%%%%%%%%%%%%%%%%%%%%%%%%%%%%%

Consider the region $R$ enclosed by the interval $q_R$ of $q^+$ 
between $e$ and $e_q^l$, pieces of these two edges and the
part $\gamma_R$ of $\gamma$ between its crossings with $e$ and $e_q^l$.
A pseudoline may intersect each of $q_R$ and $\gamma_R$ at most
once. Therefore, $q_R$ and  $\gamma_R$ are intersected 
by the same lines. This allows to construct a new curve
$\gamma'$ that starts like $\gamma$ but instead of
intersecting $q^+$ where $\gamma$ does, it 
stays on the left side of $q^+$ and only 
intersects $q^+$ at its second to last intersection, i.e.,
immediately before reaching $e_q^l$.
Figure~\ref{fig:reroute} shows an example. 

%%%%%%%%%%%%%%%%%%%%%%%%%%%%%%%%%%%%%%%%%%%%%%%%%%%%%%%%%%%
%%
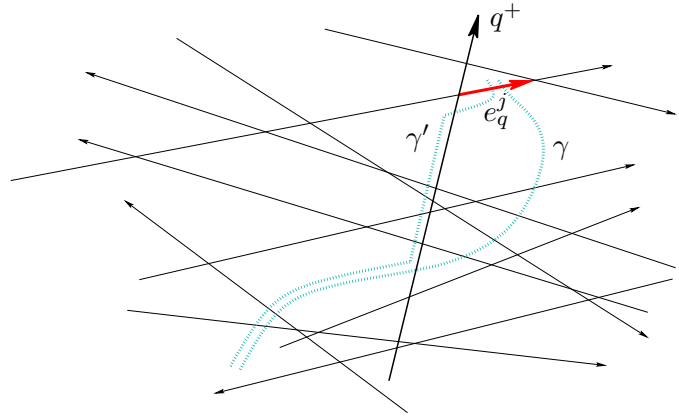
\PsFigCap{30}{reroute}{A curve $\gamma$ and the {\em normalization}
          $\gamma'$ reach $e_q^l$ with the same number of intersections.}
%%
%%%%%%%%%%%%%%%%%%%%%%%%%%%%%%%%%%%%%%%%%%%%%%%%%%%%%%%%%%%

To summarize: If $e_q^l$, the $l$th right-outgoing edge 
of $q^+$, is in $Z^{r}_k(p)$ and $l > k$,
then there is a {\em normalized} curve $\gamma'$ connecting $p$ to $e_q^l$
such that the $k$th intersection of $\gamma'$ is with $q^+$
and the $(k+1)$st intersection is with the edge $e_q^l$.
The edge $e'$ where $\gamma'$ intersects $q^+$ is 
the {\em witness} for $e_q^l \in Z^{r}_k(p)$.

\bItem If $e'$ is a witness for $e_q^l \in Z^{r}_k(p)$, then
       $e' \in Z_j(p)$ for some $j <k$.
\bItem An edge $e' \in Z_j(p)$ can be the  witness 
       for at most two edges in $Z^{r}_k(p)$. If it is 
       a witness for two, then it is a left-outgoing edge,
       hence, in $Z^{l}_j(p)$.
\medskip

\ni 
With these two observations, and collecting the $n-1$ from the beginning
of the proof we have 
$|Z^{r}_k(p)| < (n-1) + \sum_{j<k}( |Z_j(p)| + |Z^{l}_j(p)|)$.
This shows that we obtain an upper bound $|Z^{r}_k(p)| < a_k$
by solving 
$a_k = (n-1) + 3\sum_{j<k} a_j$ with the initial condition $a_1 =2(n-1)$
from the proof of the 2-D Zone Theorem.
The solution to this recurrence is $a_k = 7\cdot 4^{k-2}(n-1)$.
This is the upper bound we have stated in the lemma.
\qed

%%%%%%%%%%%%%%%%%%%%%%%%%%%%%%%%%%%%%%%%%%%%%%%%%%%%%%
\subsection{$k$-segments in arrangements}\label{ssec:k-seg}
%%%%%%%%

We now come to the proof of Theorem~\ref{thm:k-seg}.  To bound the
number of $k$-segments we proceed as follows: Consider $p$ and some
$e\in \EZone_{k-1}(p) \cup \FZone_{k-1}(p)$.  There are $n$ choices
for $p$. Theorem~\ref{thm:k-zone} bounds the number of choices for $e$.
From Lemma~\ref{lem:knuth} we get a bound on the number of
$k$-segments that have $e$ as extremal edge and the other extremal
edge on $p$.  This estimate counts every $k$-segment twice, hence
$$
\#\text{$k$-segments} <
n \cdot \big(2\cdot5\cdot 4^{k-2}(n-1)\big)\cdot k\cdot 3^{k-2}\,\frac{1}{2}
< 5k\; 12^{k-2} n^2.
$$
This shows that the constant $c$ from Theorem~\ref{thm:k-seg} can safely
be chosen as $c=12.5$.

\begin{lemma}\label{lem:knuth}
If $e$ is an edge of $\PL$ not on $p$ then then there are at most $k 3^{k-2}$
that have $e$ as extremal edge and the other extremal edge on $p$.
\end{lemma}

\Proof 
Any two curves that represent $k$-segments with $e$ and $e'$ as
extremal edges cross the same $k-2$ lines. This can be shown by
considering the region enclosed by the two curves.
It follows that we can view the $k$-segments with $e$ and $e'$ as
extremal edges as {\em cut-paths} in an arrangement of $k-2$
pseudolines. In his monograph `Axioms and Hulls' Knuth proves
the upper bound $3^n$ for the maximum number of cut-paths 
in arrangements of $n$ pseudolines, see~\cite[page 38]{k-ah-92}.

It remains to show that there are at most
$k$ choices for an edge $e_p$ on $p$ such that $e$ and~$e_p$ are the
extremal edges of a $k$-segment. 

Let $\gamma$ and $\gamma'$ be representatives of $k$-segments 
whose extremal edges are $e$ and some edges on $p$. If $\gamma$ and $\gamma'$
are disjoint there is a clear notion of which is left and which is right.
If $\gamma$ and $\gamma'$ intersect then there exist two
representatives of $k$-segments $\gamma_l$ and $\gamma_r$ such that
$\gamma_l$ is left of both and $\gamma_r$ is right of both. 
This implies that there is a leftmost~$S_l$ and a rightmost~$S_r$ 
among these $k$-segments. Consider the interval $I$
on $p$ between the edge $p \cap S_l$ and the edge $e \cap S_r$.
Every line intersecting $p$ in $I$ must have an edge on
$S_l$ or $S_r$. Since $p$ is out of question and edge $e$
belongs to both there are at most $2(k-2) +1$ vertices on $I$,
hence, $I$ contains at most $2k -2$ edges of $p$. 
Figure~\ref{fig:connecting-e} shows an example.

Finally, observe that only every other edge on $p$ can be extreme for a
$k$-segment connecting~$e$ to $p$. The only exception is with two edges
incident to the vertex where the supporting line of $e$ intersects $p$.
This shows that from the $2k -2$ candidate edges on $I$ only at most
$k$ can indeed contribute to one of the $k$-segments.
\qed

%%%%%%%%%%%%%%%%%%%%%%%%%%%%%%%%%%%%%%%%%%%%%%%%%%%%%%%%%%%
%%
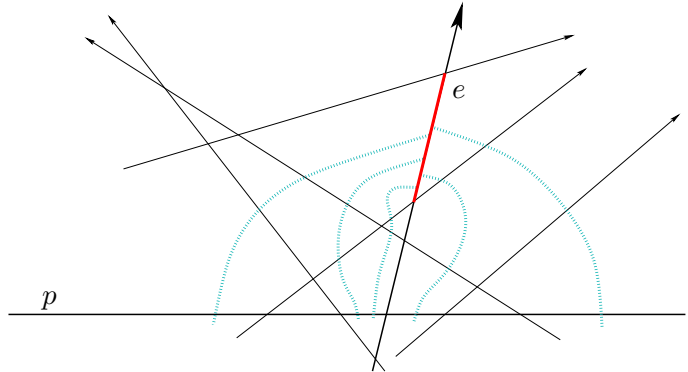
\PsFigCap{30}{connecting-e}{A line $p$, an edge $e$ and 
curves representing all 4-segments with extreme edges $e$ and and on $p$. }
%%
%%%%%%%%%%%%%%%%%%%%%%%%%%%%%%%%%%%%%%%%%%%%%%%%%%%%%%%%%%%

As mentioned by Knuth an improvement of his bound for the number of
cut-paths in an arrangement could improve the known bounds~\cite{f-nap-97}
for the number of arrangements of $n$ pseudolines.

From the point of view of our investigations the following also seems
to be very interesting: Call a set of cut-paths {\em compatible} if
they can be represented by a set of curves such that every pair of 
curves from the set has at most one intersection.

\begin{problem}
How many compatible cut-paths can an arrangement of $n$ pseudolines
have? 
\end{problem}

%%%%%%%%%%%%%%%%%%%%%%%%%%%%%%%%%%%%%%%%%%%%%%%%%%%%%%%%%%%%%%%%%%%%
\section{Chordal Graphs that are not PSI Revisited:\\
         An Application of Ramsey Theory}
%%%%%%%%%%%%%%%%%%%%%%%%%%%%%%%%%%%%%%%%%%%%%%%%%%%%%%%%%%%%%%%%%%%%

In Theorem~\ref{thm:K-3-33-not-psi} we have seen that even in
the class of chordal graphs which can be represented by
substars of a star there are graphs which are not PSI. 
If the degree of the subtrees 
is bounded by 2, then the subtrees are just paths and the graphs are PSI by
Theorem~\ref{thm:vpt-is-psi}.

What if we restrict the host tree to be a star and the substars to be
of degree three? Again we have graphs in the class that are not PSI.
This is the topic of this section.

\begin{theorem}\label{thm:indep+3-not-psi}
Let $S_n$ be the chordal graph whose
vertices are represented by all substars with three leaves and all
leaves on a star with $n$ leaves. For $n$ large enough
$S_n$ is not PSI-representable
\end{theorem}

Suppose that there is a PSI representation of $S_n$. Let $\PD$ be the
set of pseudosegments representing the leaves of the star, these segments are
pairwise disjoint. To simplify the picture we use
an homeomorphism of the plane that aligns the 
pseudosegments of $\PD$ as vertical segments of unit length which touch the 
$X$-axis with their lower endpoints at positions $1,2,\ldots,n$.
For ease of reference we will call these vertical segments {\em sticks} 
and number them such that $p_i$ is the stick containing the
point $(i,0)$.

With every ordered triple $(i,j,k)$, $1\leq i<j<k\leq n$, there is a
$3$-segment $\gamma_{ijk}$ intersecting the sticks $p_i$,
$p_j$ and $p_k$ from $\PD$. Let $\phi_{ijk}$ denote the middle of the
three sticks intersected by $\gamma_{ijk}$. We partition the
ordered triples $(i,j,k)$ into three classes depending of the position
of $\phi_{ijk}$ in the list $(p_i,p_j,p_k)$. If $\phi_{ijk}=p_i$,
i.e., the middle intersection of $\gamma_{ijk}$ is left of the other
two, we assign $(i,j,k)$ to class $[L]$. The class of $(i,j,k)$ is
$[M]$ if $\phi_{ijk}=p_j$, i.e., the middle intersection of
$\gamma_{ijk}$ is between the other two. The class of $(i,j,k)$ is
$[R]$ if $\phi_{ijk}=p_k$, i.e., the middle intersection of
$\gamma_{ijk}$ is to the right of the other two. We use this notation 
rather flexible and also write  $\gamma_{ijk} \in [X]$ or 
say that $\gamma_{ijk}$ is of class $[X]$ when the triple
$(i,j,k)$ is of class $[X]$, for $X\in L,M,R$.

%%%%%%%%%%%%%%%%%%%%%%%%%%%%%%%%%%%%%%%%%%%%%%%%%%%%%%%%%%%
%%
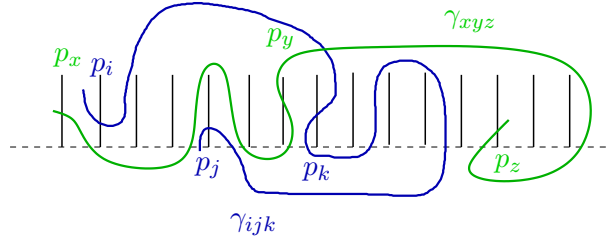
\PsFigCap{30}{lin}{Two $3$-segments 
$\gamma_{ijk}$ and $\gamma_{xyz}$. Note that $\phi_{ijk}=p_k$ and  $\phi_{xyz}=p_y$,
hence, $\gamma_{ijk}\in [R]$ and $\gamma_{xyz}\in [M]$.}
%%
%%%%%%%%%%%%%%%%%%%%%%%%%%%%%%%%%%%%%%%%%%%%%%%%%%%%%%%%%%%

Cutting $\gamma_{ijk}$ at the intersection points with the three sticks
yields two arcs $\gamma^1_{ijk}$ and $\gamma^2_{ijk}$ each connecting two sticks and
up to two ends. The ends are of no further interest. For the arcs we adopt
the convention that $\gamma^1_{ijk}$ connects $\phi_{ijk}$
to the stick further left and and $\gamma^2_{ijk}$ connects $\phi_{ijk}$
to the stick further right. In the example of Figure~\ref{fig:lin}
$\phi_{ijk} = p_k$ so that $\gamma^1_{ijk}$ is the arc connecting
$p_i$ and $p_k$ while arc $\gamma^2_{ijk}$ connects $p_j$ and $p_k$.

For a contradiction we will show that if $n$ is large enough there
have to be two 3-segments $\gamma_{ijk}$ and $\gamma_{xyz}$ such that
$\gamma^1_{ijk}$ and $\gamma^1_{xyz}$ intersect and $\gamma^2_{ijk}$
and~$\gamma^2_{xyz}$ intersect, hence, $\gamma_{ijk}$ and
$\gamma_{xyz}$ intersect at least twice which is not allowed in a
PSI-representation.

To get to that contradiction we need some control over the behavior of
3-segments between the sticks. Let $\vec{r}_x$ be a vertical ray
downwards starting at $(x,0)$, i.e., the ray pointing down from the lower end of stick
$p_x$. Let $I^s_x(ijk)$ be the number of intersections of ray $\vec{r}_x$ with
$\gamma^s_{ijk}$ and let  $J^s_x(ijk)$ be the parity of  $I^s_x(ijk)$, i.e., 
 $J^s_x(ijk) = I^s_x(ijk) \pmod 2$.  

Given an ordered 7-tuple $(a,i,b,j,c,k,d)$, 
$1\leq a < i < b < j < c < k < d\leq n$,
we call $\gamma_{ijk}$ the induced 3-segment and define 
$T^s_x = J^s_x(ijk)$. The {\em pattern}
of the tuple is the binary 8-tuple
$$
(T^1_a,T^1_b,T^1_c,T^1_d,T^2_a,T^2_b,T^2_c,T^2_d).
$$
%%%%%%%%%%%%%%%%%%%%%%%%%%%%%%%%%%%%%%%%%%%%%%%%%%%%%%%%%%%
%%
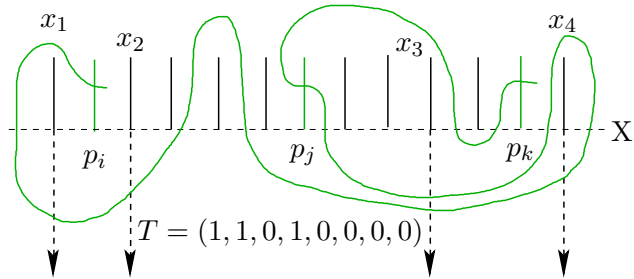
\PsFigCap{30}{ray3}{A pattern associated to a 3-segment $\gamma_{ijk}$.}
%%
%%%%%%%%%%%%%%%%%%%%%%%%%%%%%%%%%%%%%%%%%%%%%%%%%%%%%%%%%%%

The {\em color} of a 7-tuple $(a,i,b,j,c,k,d)$ is the pair consisting
of the class of the induced 3-segment and the pattern. The 7-tuples
are thus colored with the 768 colors from the set $[3]\times {\bf
  2}^8$. Ordered 7-tuples and 7-element subsets of $[n]$ are
essentially the same. Therefore we can apply the hypergraph Ramsey
theorem with parameters $768,7,13$.

\begin{theorem}[Hypergraph Ramsey]
For every choice of 
numbers $r,p,k$ there exists a number $n$ such that whenever
$X$ is an $n$-element set and $c$ is a coloring of the
system of all $p$-element subsets of $X$ using $r$ colors, i.e. $c: 
{X \choose p} \rightarrow \{1,2,..,r\}$, then there is an $k$-element
subset $Y\subseteq X$ such that all the $p$-subsets in $Y\choose p$
have the same color.
\end{theorem}

Given two curves $\gamma$ and $\gamma'$, closed or not, we let
$X(\gamma,\gamma')$ be the number of crossing points of 
the two curves. The following is essential for the argument:

\begin{fact}\label{fact:A}
If $\gamma$ and $\gamma'$ are closed curves, then
$X(\gamma,\gamma') \equiv 0 \pmod 2$.
\end{fact}

With an arc $\gamma_{ij}$ connecting sticks $p_i$ and $p_j$ we
associate a closed curve $\bga_{ij}$ as follows: At the intersection
of $\gamma_{ij}$ with either of the sticks we append long vertical
segments and connect the lower endpoints of these two segments
horizontally. The union of the three connecting segments will be
called the {\em bow} $\beta_{ij}$ of the curve $\bga_{ij}$. If this
construction is applied to several arcs we assume that the vertical
segments of the bows are long enough as to avoid any intersection between
the arcs and the horizontal part of the bows.

Given $\bga_{ij}$ and $\bga_{xy}$ we can count their crossings in parts:
$$
X(\bga_{ij},\bga_{xy}) = X(\gamma_{ij},\gamma_{xy}) + X(\gamma_{ij},\beta_{xy}) +
                         X(\beta_{ij},\gamma_{xy})  + X(\beta_{ij},\beta_{xy})
$$
With Fact~\ref{fact:A} we obtain
\begin{fact}\label{fact:B}
$X(\gamma_{ij},\gamma_{xy}) \equiv
 X(\gamma_{ij},\beta_{xy}) + X(\beta_{ij},\gamma_{xy})  + X(\beta_{ij},\beta_{xy})
 \pmod 2.$
\end{fact}

The application of the Ramsey theorem left us with a uniformized configuration.
We have kept only a subset $Y$ of sticks such that all  3-segments connecting three
of them are of the same class and all 7-tuples on $Y$ have the same
pattern $T= (T^1_1,T^1_2,T^1_3,T^1_4,T^2_1,T^2_2,T^2_3,T^2_4)$.

The uniformity allows to apply Fact~\ref{fact:B} to detect 
intersections of arcs of type $\gamma^s_{ijk}$.
Depending on the entries of the pattern $T$ we choose two appropriate
3-segments $\gamma_{ijk}$ and $\gamma_{xyz}$ and show that they
intersect twice. Assume that the class of all 3-segments is $[L]$ or$[M]$,
hence there is an arc connecting the two sticks with smaller indices.
Let $\gamma_{ij} = \gamma^1_{ijk}$, and $\gamma_{xy} = \gamma^1_{xyz}$.

\begin{lemma}\label{lem:alt}
If $T^1_1 = T^1_3$ and $i < x < j < y < k$, then there is an intersection between
the arcs $\gamma_{ij}$ and $\gamma_{xy}$.
\end{lemma}

\Proof
We evaluate the right side of the congruence given in Fact~\ref{fact:B}.

$X(\gamma_{ij},\beta_{xy})$ is the number of intersections of
arc $\gamma_{ij}$ with the bow connecting $p_x$ and $p_y$.
These intersections happen on the vertical part, hence on the
rays $\vec{r}_x$ and $\vec{r}_y$. The parity of these intersections 
can be read from the pattern. The position of $x$ between $i$ and $j$
implies $T_x^1 = T^1_2$ and the position of $y$ between $j$ and $k$ 
implies $T_y^1 = T^1_3$. Hence, $X(\gamma_{ij},\beta_{xy})\equiv
T^1_2 + T^1_3 \pmod 2$.

From the positions of $i$ left of $x$ and of $j$ between $x$ and $y$ we
conclude that $X(\beta_{ij},\gamma_{xy}) \equiv T^1_1 + T^1_2 \pmod 2$.

Since the pairs $ij$ and $xy$ interleave the two bows are intersecting, i.e., 
$X(\beta_{ij},\beta_{xy}) = 1$.

Together this yields
$X(\gamma_{ij},\gamma_{xy}) \equiv T^1_2 + T^1_3 + T^1_1 + T^1_2 + 1 \pmod 2$.
With  $T^1_1 = T^1_3$ we see that  $X(\gamma_{ij},\gamma_{xy})$ is odd,
hence, there is at least one intersection between the arcs.
\qed

\begin{lemma}\label{lem:non-alt}
If $T^1_1 \neq T^1_3$ and $x < i < j < y < k$, then there is an intersection between
the arcs $\gamma_{ij}$ and $\gamma_{xy}$.
\end{lemma}

\Proof
Since $x$ is left of $i$ and $y$ is between $j$ and $k$ we obtain
$X(\gamma_{ij},\beta_{xy}) \equiv T^1_1 + T^1_3 \pmod 2$.
Both $i$ and $j$ are between $x$ and $y$, thus 
$X(\beta_{ij},\gamma_{xy}) \equiv T^1_2 + T^1_2 \equiv 0 \pmod 2$.
For the bows we observe that either they don't intersect or they intersect
twice, in both cases $X(\beta_{ij},\beta_{xy}) \equiv 0 \pmod 2$.

Put together $X(\gamma_{ij},\gamma_{xy}) \equiv T^1_1 + T^1_3 \pmod 2$.
With  $T^1_1 \neq T^1_3$ we see that  $X(\gamma_{ij},\gamma_{xy})$ is odd,
hence, there is at least one intersection between the arcs.
\qed

Now consider the case where the class of all 3-segments is $[M]$.
In addition to the arcs $\gamma_{ij}$ and $\gamma_{xy}$ we have the arcs
$\gamma_{jk} = \gamma^2_{ijk}$, and $\gamma_{yz} = \gamma^2_{xyz}$.
The following two lemmas are counterparts to lemmas~\ref{lem:alt} and~\ref{lem:non-alt}
they show that depending on the parity of $T^2_1+T^2_3$ an alternating or a non-alternating
choice of $jk$ and $yz$ force an intersection of the arcs $\gamma_{jk}$ and $\gamma_{yz}$.
For the proofs note that reflection at the $y$-axis keeps class $[M]$ invariant
but exchanges the first and the second arc, the relevant effect on the pattern is
$T^1_1 \leftrightarrow T^2_4$ and $T^1_3 \leftrightarrow T^2_2$. 

\begin{lemma}\label{lem:alt2}
If $T^2_2 = T^2_4$ and $x < j < y < k < z$, then there is an intersection between
the arcs $\gamma_{jk}$ and $\gamma_{yz}$.
\end{lemma}

\begin{lemma}\label{lem:non-alt2}
If $T^2_2 \neq T^2_4$ and $x < j < y < z < k$, then there is an intersection between
the arcs $\gamma_{jk}$ and $\gamma_{yz}$.
\end{lemma}

The table below shows that it is possible to select $ijk$ and $xyz$ out of six
numbers such that the positions of $ij$ and $xy$ respectively $jk$ and $yz$ are
any combination of alternating and non-alternating. Hence, according to the lemmas
we have at least two intersections between 3-segments $\gamma_{ijk}$ and $\gamma_{xyz}$
chosen appropriately depending on the entries of pattern $T$. We
represent elements of $ijk$ by a box $\Box$ and elements of $xyz$ by circles $\bullet$.

$$
\begin{tabular}{>{$}c<{$\kern-7pt}>{$}c<{$\kern-7pt}>{$}c<{$\kern-7pt}>{$}c<{$\kern-7pt}>{$}c<{$\kern-7pt}>{$}c<{$\kern-7pt}>{\hspace{6mm}}ll}
\Box  &  \bullet  &  \Box  &  \bullet  &  \Box  &  \bullet  &
   alt / alt          &[$T^1_1 = T^1_3$ and $T^2_2 = T^2_4$]\\
\Box  &  \bullet  &  \Box  &  \bullet  &  \bullet  &  \Box  &
   alt / non-alt      &[$T^1_1 = T^1_3$ and $T^2_2 \neq T^2_4$]\\
\bullet  &  \Box  &  \Box  &  \bullet  &  \Box  &  \bullet  &
   non-alt / alt      &[$T^1_1 \neq T^1_3$ and $T^2_2 = T^2_4$]\\
\bullet  &  \Box  &  \Box  &  \bullet  &  \bullet  &  \Box  &
   non-alt / non-alt  &[$T^1_1 \neq T^1_3$ and $T^2_2 \neq T^2_4$]\\
\end{tabular}
$$

Now consider the case where the class of all 3-segments is $[L]$.  In
addition to the arcs $\gamma_{ij}$ and $\gamma_{xy}$ we have the arcs
$\gamma_{ik} = \gamma^2_{ijk}$, and $\gamma_{xz} = \gamma^2_{xyz}$.
The following two lemmas show that depending on the parity of
$T^2_1+T^2_3+T^2_3+T^2_4$ an alternating or a non-alternating choice
of $ik$ and $xz$ force an intersection of the arcs $\gamma_{ik}$ and
$\gamma_{xz}$.

\begin{lemma}\label{lem:alt3}
If $T^2_1+T^2_3+T^2_3+T^2_4 \equiv 0 \pmod 2$ and $i < x < \{j,y\} < k < z$, 
then there is an intersection between
the arcs $\gamma_{ik}$ and $\gamma_{xz}$.
\end{lemma}

\Proof
Since $x$ is between $i$ and $j$ and $z$ is right of $k$ we obtain
$X(\gamma_{ik},\beta_{xz}) \equiv T^2_2 + T^2_4 \pmod 2$.
Since $i$ is left of $x$ and $k$ is between $y$ and $z$ we obtain
$X(\beta_{ik},\gamma_{xz}) \equiv T^2_1 + T^2_3 \pmod 2$.
Since the pairs $ik$ and $xz$ interleave the two bows are intersecting, i.e., 
$X(\beta_{ik},\beta_{xz}) = 1$.

Put together $X(\gamma_{ij},\gamma_{xy}) \equiv T^2_1+T^2_3+T^2_3+T^2_4 + 1 \pmod 2$.
Hence there is at least one intersection between the arcs.
\qed

\begin{lemma}\label{lem:alt4}
If $T^2_1+T^2_3+T^2_3+T^2_4 \equiv 1 \pmod 2$ and $i < x < \{j,y\} < z < k$, 
then there is an intersection between
the arcs $\gamma_{ik}$ and $\gamma_{xz}$.
\end{lemma}

\Proof
Since $x$ is between $i$ and $j$ and $z$ is between $j$ and $k$ we obtain
$X(\gamma_{ik},\beta_{xz}) \equiv T^2_2 + T^2_3 \pmod 2$.
Since $i$ is left of $x$ and $k$ is right of $z$ we obtain
$X(\beta_{ik},\gamma_{xz}) \equiv T^2_1 + T^2_4 \pmod 2$.
For the bows we observe that either they don't intersect or they intersect
twice, hence $X(\beta_{ik},\beta_{xz}) \equiv 0 \pmod 2$.

Put together $X(\gamma_{ij},\gamma_{xy}) \equiv T^2_1+T^2_3+T^2_3+T^2_4 \pmod 2$.
Hence there is at least one intersection between the arcs.
\qed

As in the previous case we provide a table 
showing that it is possible to select $ijk$ and $xyz$ out of six
numbers such that the positions of $ij$ and $xy$ respectively $ik$ and $xz$ are
any combination of alternating and non-alternating. We
represent elements of $ijk$ by a box $\Box$ and elements of $xyz$ by circles $\bullet$.

$$
\begin{tabular}{>{$}c<{$\kern-7pt}>{$}c<{$\kern-7pt}>{$}c<{$\kern-7pt}>{$}c<{$\kern-7pt}>{$}c<{$\kern-7pt}>{$}c<{$\kern-7pt}>{\hspace{4mm}}ll}
\Box  &  \bullet  &  \Box  &  \bullet  &  \Box  &  \bullet  &
   alt / alt          &[$T^1_1 = T^1_3$ and $T^2_1+T^2_3+T^2_3+T^2_4 \equiv 0$]\\
\Box  &  \bullet  &  \Box  &  \bullet  &  \bullet  &  \Box  &
   alt / non-alt      &[$T^1_1 = T^1_3$ and $T^2_1+T^2_3+T^2_3+T^2_4 \equiv 1$]\\
\bullet  &  \Box  &  \Box  &  \bullet  &  \bullet  &  \Box  & 
   non-alt / alt      &[$T^1_1 \neq T^1_3$ and $T^2_1+T^2_3+T^2_3+T^2_4 \equiv 0$]\\
\bullet  &  \Box  &  \Box  &  \bullet  &  \Box  &  \bullet  &
   non-alt / non-alt  &[$T^1_1 \neq T^1_3$ and $T^2_1+T^2_3+T^2_3+T^2_4 \equiv 1$]\\
\end{tabular}
$$

To deal with the case where the class of all 3-segments is $[R]$ we
refer to symmetry. Reflecting the picture at the $y$-axis yields a
configuration which is in class $[L]$. Hence it is impossible to have
a uniform configuration on six or more sticks. Well, the definition of
the pattern alone involves seven sticks. This is why we said that we
want to have a uniform family on 13 sticks, if we use the odd numbered
sticks for the choice of triples there are enough candidates to
complete the selection of seven positions for a pattern. This
completes the proof of Theorem~\ref{thm:indep+3-not-psi}.

%%%%%%%%%%%%%%%%%%%%%%%%%%%%%%%%%%%%%%%%%%%%%%%%%%%%%%%%%%%%%%%%%%%%

\section{Conclusion}

Besides the problems that have already been presented in the text
we would like to mention some more.
In a large part of the paper we been considering questions 
of the following general type: How big a
family $W$ of 3-segments can be that live on a base arrangement $B$
of size $n$?

\bItem
In Theorem~\ref{thm:K-3-33-not-psi} we requested that the 3-segments 
are disjoint.

\bItem
In Theorem~\ref{thm:k-seg} we requested that $B$ is given by an
arrangement of pseudolines.

\bItem
In Theorem~\ref{thm:indep+3-not-psi} we requested that the segments 
in $B$ are disjoint and that the segments in~$W$ are compatible.
\medskip

\noindent
The theorems give upper bounds. Good lower bound constructions might
give some additional insight that can help improve the upper bounds.
In fact we think that in the situation of
Theorem~\ref{thm:indep+3-not-psi} the true size of $W$ should be in
$o(n^{2+\epsilon})$.

\small
\global\advance\baselineskip 2pt
% *************************************************************
%%\bibliography{schnyderwoods}
\bibliographystyle{my-siam}
\bibliography{psi}

\end{document}